\documentclass[12pt]{article}
\usepackage{amsxtra,amssymb,amsthm,amsmath,latexsym}

\theoremstyle{plain}
\newtheorem{theorem}{Theorem}

\newtheorem{remark}{Remark}
\newtheorem*{A}{Assumption A}

\numberwithin{equation}{section}

\newcommand{\refT}[1]{Theorem~\ref{T:#1}}
\newcommand{\refS}[1]{Section~\ref{S:#1}}

\def\qed{{\hfill $\Box$}}

\def\dotg{\dot{g}}

\def\dotw{\dot{w}}

\def\dotw{\dot{w}}
\def\ve{{\varepsilon}}

\def\oH1{{\overset{\circ}{H}\kern-.04in{}^1}}

\def\loc{{\,loc\,}}
\def\bee{\begin{equation*}}
\def\eee{\end{equation*}}
\def\be{\begin{equation}}
\def\ee{\end{equation}}

\begin{document}

\title{Dynamical systems method (DSM) for nonlinear equations in Banach 
spaces}

\author{
A.G. Ramm\\
 Mathematics Department, Kansas State University, \\
 Manhattan, KS 66506-2602, USA\\
ramm@math.ksu.edu\\
http://www.math.ksu.edu/\,$\widetilde{\ }$\,ramm}
\date{}
\maketitle\thispagestyle{empty}

\begin{abstract}
\footnote{Math subject classification: 47J05, 47J06, 47J25 }
\footnote{key words: dynamical systems method, operator equations, 
ill-posed problems, nonlinear problems  }
Let $F:X\to X$ be a $C^2_\loc$ map in a Banach space $X$, and $A$ be its
Fr\`echet derivative at the element $w:=w_\ve$, which solves the problem
$(\ast)\ \dotw=-A^{-1}_\ve(F(w)+\ve w)$,
$w(0)=w_0$, where $A_\ve:=A+\ve I$.
Assume that $\|A^{-1}_\ve\|\leq c \ve^{-k}$, $0<k\leq 1$, $0<\ve>\ve_0$.
Then $(\ast)$ has a
unique global solution, $w(t)$, there exists $w(\infty)$, and
$(\ast\ast)\ F(w(\infty))+\ve w(\infty)=0$.
Thus the DSM (Dynamical Systems Method) is justified for equation 
$(\ast\ast)$.
The limit of $w_\ve$ as $\ve\to 0$ is studied.
\end{abstract}

\section{Introduction}\label{S:1}
In \cite{R456}  the DSM (Dynamical Systems Method) was developed for solving
operator equations 
\be\label{e1.1} F(u)=0\ee
in a Hilbert space. In this note we generalize the DSM for Banach spaces and a
more general class of nonlinear operators.

Let $X$ be a Banach space, {\it not necessarily reflexive}, and $F:X\to X$ 
be a
$C^2_\loc$ map. This means that $F$ is twice Fr\`echlet differentiable and

\be\label{e1.2}
\sup_{u\in B(u_0,R)} \|F^{(j)}(u)\|\leq M_j(R) ,\quad j=1,2,\ee
where $B(u_0,R):=\{u:\|u-u_0\|\leq R\}$. Let $F'(u):=A:=A(u), A_\ve:=A+\ve I$,
where $I$ is the identity operator.

\begin{A}
\be\label{e1.3} 
\|A^{-1}_\ve\|\leq c_0\ve^{-k}, \qquad \ve\in(0,\ve_0), \qquad 
k=const>0,\ee
where $c_0=const>0$, $\ve_0>0$ is an arbitrary small fixed number.
\end{A}

\begin{theorem}\label{T:1}
If \eqref{e1.2} and \eqref{e1.3} hold, then equation
\be\label{e1.4}
F(u)+\ve u=0, \qquad \ve\in(0,\ve_0),\ee
has a solution.
\end{theorem}

In \refS{2} this result is proved.

In Theorem 2 of \refS{3} conditions for the convergence $u:= u_\ve\to u_0$ 
as $\ve\to\infty$ are given, where $u_0$ solves \eqref{e1.1}.

\section{Proof of \refT{1}}\label{S:2} Consider the equation
\be\label{e2.1} \dotw=-A^{-1}_\ve (w) [F(w)+\ve w], \qquad w(0)=w_0,
\quad \ve\in(0,\ve_0),
\ee
where $\dotw=\frac{dw}{dt}$ is the strong derivative, and $w_0\in X$ 
is
arbitrary.  Let $h\in X^\ast$ be an arbitrary linear bounded functional on
$X$. Define $g(t):=(F(w)+\ve w,h)$, where $(u,h)$ is the value of the
functional $h$ on the element $u\in X$, and $w=w(t)$ is the local solution
to \eqref{e2.1}. From the assumptions \eqref{e1.2} and \eqref{e1.3} it
follows that the right-hand side of \eqref{e2.1} is locally Lipschitz, so
\eqref{e2.1} has a unique local solution $w$. We wish to justify the DSM
for solving equation \eqref{e1.4}. 

The DSM consists of: 

a) proving that
$w(t)$ exists globally, i.e., $\forall t>0$, 

b) the limit
$\lim_{t\to\infty}w(t):=w(\infty)$ exists, 

and 

c) $w(\infty)$ solves
\eqref{e1.4}.

To prove a), b) and c), we start with the equation
$\dotg=(A_\ve(w)\dotw,h)=-g$, which implies:
\be\label{e2.2} g(t)=g(0)e^{-t}, \qquad
\|g(t)\|=\|g(0)\|e^{-t}.\ee Thus \be\label{e2.3} \|F(w)+\ve w\|
=\sup_{\|h\|\leq 1} |g(t)| \leq \|F(w_0)+\ve w_0\|e^{-t} :=F_0 e^{-t}, \ee
and \eqref{e2.1} implies: \be\label{e2.4} \|\dotw\| \leq c_0\ve^{-k} F_0
e^{-t}.\ee From \eqref{e2.4} it follows that
$\|\dotw\|\in L^1(0,\infty)$. This and the Cauchy test for the existence
of the limit $w(\infty):=\lim_{t\to \infty}w(t)$ imply 
that
$w(\infty)$ exists,  \be\label{e2.5} \|w(t)-w(\infty)\|\leq
F_0c_0\ve^{-k} e^{-t},\ee and \be\label{e2.6} \|w(t)-w_0\| \leq F_0
c_0\ve^{-k}.\ee From \eqref{e2.4} and \eqref{e2.1}, passing to the limit
$t\to\infty$, one gets \bee 0=-A^{-1}_\ve (w(\infty)) [F(w(\infty))+\ve
w(\infty)].\eee This implies that $w(\infty)$ solves \eqref{e1.4}. The DSM
is justified. \refT{1} is proved. \qed

\begin{remark}\label{R:1} The solution $w(\infty)$ depends on the choice
of $w_0$. Equation \eqref{e1.4} may have many solutions under the
assumptions of \refT{1}. \end{remark}

\section{Limiting behavior of the solution.}\label{S:3} Denote
$w(\infty):=v_\ve:=v$. We want to give conditions sufficient for the
existence of the limit \be\label{e3.1} \lim_{\ve\to 0} v_\ve=y,\ee where
$y$ solves equation \eqref{e1.1}.

First, note that equation \eqref{e1.4} can be solvable for any
$\ve\in(0,\ve_0)$, but equation \eqref{e1.1} may have no solution. For
example, let $F(u)=Au-f$, where $A\geq 0$ is a bounded selfadjoint
operator in $X=H$, where $H$ is a Hilbert space and $f\not\in R(A)$, where
$R(A)$ is the range of $A$. The equation $A_\ve v_\ve-f=0$ has a unique
solution for any $\ve>0$, but the limiting equation $Ay-f=0$ has no
solution.
That is why we assume that equation \eqref{e1.1} is solvable: $F(y)=0$.
If $F(y)=0$,
then \bee 0=F(v)+\ve v-F(y)=F'(y)(v-y)+R+\ve(v-y)+\ve y.\eee 
By Taylor's formula one has $F(v)-F(y)=F'(y)(v-y)+R$. Let $z:=v-y$. Then
\be\label{e3.2} A_\ve z+R=-\ve y, \qquad A:=F'(y).
\ee Assume that \be\label{e3.3} y=A\psi,
\qquad \|\psi\|\ll 1,\ee where $\|\psi\|\ll 1$ means that
$\|\psi\|$ is sufficiently small (see \eqref{e3.8} below).

Then \eqref{e3.2} is equivalent to \be\label{e3.4} z=-A^{-1}_\ve R-\ve
A^{-1}_\ve A\psi:=T(z),\ee where \be\label{e3.5}
R:=F(v)-F(y)-F'(y)z=\int^1_0 ds(1-s) F''(y+sz)zz.\ee Let us check that the
map $T$ maps a ball $B(0,z):=B_r:=\{u:\|z\|\leq r\}$, $z=v-y$, into itself
and is a contraction in $B_r$ for a suitable $r>0$. Indeed,
\be\label{e3.6} \|T(z)\|\leq \frac{c_0}{\ve^k} \ \frac{M_2}{2}
r^2+\ve\|\psi\|\leq r,\ee provided that \be\label{e3.7} r=\frac{\ve^k}{c_0
M_2} \left(1-\sqrt{1-2c_0M_2\|\psi\| \ve^{1-k}}\right), \ee and
\be\label{e3.8} \rho:=2c_0M_2\|\psi\| \ve^{1-k}<1.\ee Condition 
\eqref{e3.8}
is satisfied if $k<1$ and $\ve$ is sufficiently small, or if $k=1$ and
$\|\psi\|$ is sufficiently small. If $k>0$ then $r=r(\ve)\to 0$ as $\ve\to
0$, and $T$ maps $B_{r(\ve)}$ into itself.

Let us check the contraction mapping property. Let $z,p\in B_r$. Then,
using \eqref{e3.5}, one gets \be\label{e3.9} \begin{aligned} \|T(z)-T(p)\|
&\leq {\frac{c_0}{\ve^k}} \|R(z)-R(p)\| \\ & \leq {\frac{c_0}{\ve^k}}
   \int^1_0 ds(1-s) \left[ \|F''(y+sz)-F''(y+sp)\| r^2 +2M_2 r \|z-p\|\right]\\
&\leq {\frac{c_0}{\ve^k}}
    \int^1_0 ds(1-s) \left[ s M_3 r^2 + 2M_2 r\right]\|z-p\| \\
&\leq {\frac{c_0}{\ve^k}} \left(\frac{M_3}{6} r^2+M_2 r\right) \|z-p\|. 
\end{aligned}\ee
Thus $T$ is a contraction on $B_r$ if 
\be\label{e3.10}
\eta:=\frac{c_0}{\ve^k} \left(\frac{M_3r^2}{6} + M_2r\right)\leq q<1.\ee

If \eqref{e3.7} and \eqref{e3.8} hold, then 
\be\label{e3.11}
\eta=O(\ve^k)+1-\rho:=q<1 \ee
if $\ve\in(0,\ve_0)$ is sufficiently small.

We have proved:

\begin{theorem}\label{T:2} Assume that: 1) equation \eqref{e1.1} is
solvable, $F(y)=0$, 2) \eqref{e1.2} holds for $j\leq 3$, 3)\eqref{e1.3}
holds, 4) \eqref{e3.3} holds and 5) \eqref{e3.8} holds.  Then there exists
and is unique a solution $v_\ve$ to equation \eqref{e1.4} such that
\be\label{e3.12} \|v_\ve-y\|=O(\ve^k).\ee \end{theorem}

\begin{remark}\label{R:2} One may drop assumption \eqref{e3.3} and
consider in place of equation \eqref{e1.4} the following one:
\be\label{e3.13} F(p)+\ve(p-q)=0, \qquad p:=p_\ve, \ee where $q$ is any
element such that \be\label{e3.14} y-q=A\psi, \qquad \|\psi\|\ll
1.\ee Then \eqref{e3.12} holds with $p_\ve$ in place of $v_\ve$ and the
proof is essentially the same. Note that \eqref{e3.14} holds, with
$A:=F'(y)$ if $AB_r\cap (B_a\backslash\{0\})\not=\emptyset$ for any
$r\in(0,r_0)$, where $r_0>0$ is a fixed number. \end{remark}

\end{document}